\documentclass[a4, 11pt, leqno]{article}
\usepackage{amsmath}
\usepackage{amssymb}
\newtheorem{df}{Definition}[section]

\newtheorem{Th}[df]{Theorem}

\begin{document}

\title{\bf Motion of a Vortex Filament in an External Flow}
\author{Masashi A{\sc iki} and Tatsuo I{\sc guchi}}
\date{}

\maketitle
\vspace*{-0.5cm}

\begin{abstract}
We consider a nonlinear model equation describing the motion of a vortex filament
immersed in an incompressible and inviscid fluid. In the present problem setting, we 
also take into account the effect of external flow. We prove the unique solvability,
locally in time, of an initial value problem posed on the one dimensional torus. The  
problem describes the motion of a closed vortex filament.
\end{abstract}

\section{Introduction and Problem Setting}
A vortex filament is a space curve on which the vorticity of the fluid is concentrated. 
Vortex filaments are used to model very thin vortex structures such as vortices that 
trail off airplane wings or propellers. 
In this paper, we prove the solvability of the following initial value problem
which describes the motion of a closed vortex filament.
\begin{eqnarray}
\left\{
\begin{array}{ll}
\displaystyle \mbox{\mathversion{bold}$x$}_{t} = \frac{\mbox{\mathversion{bold}$x$}_{\xi }\times 
\mbox{\mathversion{bold}$x$}_{\xi \xi}}{|\mbox{\mathversion{bold}$x$}_{\xi}|^{3}}
+ \mbox{\mathversion{bold}$F$}(\mbox{\mathversion{bold}$x$},t), &
\xi \in \mathbf{T}, \ t>0, \\[2ex]
\mbox{\mathversion{bold}$x$}(\xi ,0) = \mbox{\mathversion{bold}$x$}_{0}(\xi), & \xi \in \mathbf{T},
\end{array}\right.
\label{ext}
\end{eqnarray}
where \( \mbox{\mathversion{bold}$x$}(\xi, t)= (x_{1}(\xi ,t), x_{2}(\xi ,t), x_{3}(\xi ,t))\) is the 
position of the vortex filament parametrized by \( \xi \) at time \( t\), 
the symbol \( \times \) is the 
exterior product in the three dimensional Euclidean space, \( \mbox{\mathversion{bold}$F$}(\cdot ,t) \) is a
given external flow field, \( \mathbf{T}\) is the one dimensional torus \( \mathbf{R}/\mathbf{Z}\), 
and subscripts are differentiations with the respective variables.
Problem (\ref{ext}) describes the motion of a closed vortex filament under the influence of 
external flow. Such a setting can be seen as an idealization of the motion of a bubble ring in water, where 
the thickness of the ring is taken to be zero and some environmental flow is also present. 
Many other phenomena can be modeled by a vortex ring or a closed vortex filament and are
important in both application and theory. Here, we make the distinction between a vortex ring and a 
closed vortex filament. A vortex ring is a closed vortex tube, in the shape of a torus, which has a finite core thickness. A closed vortex filament is a
closed curve, which can be regarded as a vortex ring with zero core thickness. 

The equation in problem (\ref{ext}) is a generalization of an equation called the 
Localized Induction Equation (LIE) given by 
\begin{eqnarray*}
\mbox{\mathversion{bold}$x$}_{t} = \mbox{\mathversion{bold}$x$}_{s}\times \mbox{\mathversion{bold}$x$}_{ss},
\end{eqnarray*}
which is derived by applying the so-called localized induction approximation to the Biot--Savart integral.
Here, \( s\) is the arc length parameter of the filament. The LIE was first derived by 
Da Rios \cite{31} and was re-derived twice independently by Murakami et al. \cite{33} and by 
Arms and Hama \cite{32}. Many researches have been done on the LIE and many results have been obtained.
Nishiyama and Tani \cite{5,13} proved the unique solvability of the initial value problem in Sobolev 
spaces. Koiso considered a geometrically generalized setting in which he rigorously proved the 
equivalence of the LIE and a nonlinear Schr\"odinger equation. This equivalence was first shown by 
Hasimoto \cite{14} in which he studied the formation of solitons on a vortex filament. He defined a transformation of variable known as the Hasimoto transformation to transform the LIE into a
nonlinear Schr\"odinger equation. 
The Hasimoto transformation was proposed by Hasimoto \cite{14}
and is a change of variable given by 
\[ \psi = \kappa \exp \left( {\rm i} \int ^{s}_{0} \tau \,{\rm d}s \right), \]
where \( \kappa \) is the curvature and \( \tau \) is the torsion of the filament.
Defined as such, it is well known that \( \psi \) satisfies the 
nonlinear Schr\"odinger equation given by 
\begin{eqnarray}
{\rm i}\frac{\partial \psi}{\partial t}
= \frac{ \partial ^{2}\psi }{\partial s^{2}} + \frac{1}{2} \left| \psi \right| ^{2}\psi .
\label{NLS}
\end{eqnarray}
The original transformation proposed by Hasimoto uses the torsion of the filament in its definition,
which means that the transformation is undefined at points where the curvature of the filament is zero. 
Koiso \cite{12} constructed a transformation, sometimes referred to as the 
generalized Hasimoto transformation, and gave a mathematically rigorous 
proof of the equivalence of the LIE and (\ref{NLS}). More recently, 
Banica and Vega \cite{15,16} and Guti\'errez, Rivas, and Vega \cite{17} constructed and analyzed
a family of self-similar solutions of the LIE which forms a corner in finite time. 
The authors \cite{11} proved the unique solvability of an initial-boundary value problem for the LIE
in which the filament moved in the three-dimensional half space. Nishiyama and Tani \cite{5} also
considered initial-boundary value problems with different boundary conditions. These results 
fully utilize the property that a vortex filament moving according to the LIE doesn't stretch 
and preserves its arc length parameter. 
This is not the case when we take into account the effect of external flow.
We mention here that there is another model 
describing the motion of singular vortices
such as a vortex filament, which is derived by
approximating the Biot--Savart integral. The model is called the 
Rosenhead model, which was proposed by 
Rosenhead \cite{18}, and is derived by 
desingularizing the Biot--Savart integral.
The time-local and time-global unique solvability of the 
Rosenhead model for a closed vortex filament is obtained in 
Berselli and Bessaih \cite{20} and 
Berselli and Gubinelli \cite{19} under the assumption that 
there is no external flow.

The LIE can be naturally generalized to take into account the effect of external flow.
The model equation is given by 
\begin{eqnarray}
\mbox{\mathversion{bold}$x$}_{t} = \frac{\mbox{\mathversion{bold}$x$}_{\xi }\times 
\mbox{\mathversion{bold}$x$}_{\xi \xi}}{|\mbox{\mathversion{bold}$x$}_{\xi}|^{3}}
+ \mbox{\mathversion{bold}$F$}(\mbox{\mathversion{bold}$x$},t).
\label{gLIE}
\end{eqnarray}
Here, the parametrization of the filament has been changed to \( \xi \) because, unlike 
the LIE, a vortex filament moving according to (\ref{gLIE}) stretches in general and the 
arc length is no longer preserved. It is worth mentioning that if the Jacobi matrix of
\( \mbox{\mathversion{bold}$F$}\) is skew-symmetric, which amounts to assuming that the 
effect of external flow consists only of translation and rigid body rotation, then the 
solvability for (\ref{gLIE}) can be considered in the 
same way as for the LIE. This is because if the Jacobi matrix is skew-symmetric, then the filament
no longer can stretch, and the techniques used in the analysis of the LIE can be utilized for 
(\ref{gLIE}). Thus, in what follows, we do not assume any structural conditions on 
\( \mbox{\mathversion{bold}$F$} \). 

Regarding the solvability of (\ref{gLIE}), Nishiyama \cite{10} proved the existence of 
weak solutions to initial and initial-boundary value problems in Sobolev spaces. The solutions obtained by 
Nishiyama are weak in the sense that the uniqueness of the solution is not known, but the 
equation is satisfied in the point wise sense almost everywhere. The result presented in this 
paper is an extension of Nishiyama's result for the initial value problem, and we proved the 
unique solvability in higher order Sobolev spaces. 

The contents of the rest of the paper are as follows. In Section 2, we introduce notations used 
in this paper and state our main theorem. In Section 3, we give a description for the 
construction of the solution. In Section 4, we obtain energy estimates of the solution in Sobolev spaces. 
The derivation of the energy estimate is the most crucial part of the proof of the main theorem. 
In Section 5, we prove the uniqueness of the solution along the line of the energy estimate 
carried out in Section 4. 
Finally in Section 6, we give concluding remarks.


\section{Function Spaces, Notations, and Main Theorem}

We define some function spaces that will be used throughout this paper, and 
introduce notations associated with the spaces.
For a non-negative integer \( m\), and \( 1\leq p \leq \infty \), \( W^{m,p}(\mathbf{T})\) is the Sobolev space 
containing all real-valued functions that have derivatives in the sense of distribution up to order \( m\) 
belonging to \( L^{p}(\mathbf{T})\).
We set \( H^{m}(\mathbf{T}) := W^{m,2}(\mathbf{T}) \) as the Sobolev space equipped with the usual inner product.
The norm in \( H^{m}(\mathbf{T}) \) is denoted by \( \| \cdot \|_{m} \) and we simply write \( \| \cdot \| \) 
for \( \|\cdot \|_{0} \). Otherwise, for a Banach space \( X\), the norm in \( X\) is written as \( \| \cdot \| _{X}\).
The inner product in \( L^{2}(\mathbf{T})\) is denoted by \( (\cdot ,\cdot )\).
We also make use of the Fourier series expansion for functions belonging to \( H^{m}(\mathbf{T})\).
For \( u\in H^{m}(\mathbf{T})\) and \( k\in \mathbf{Z}\), we define the 
\( k\)-th Fourier coefficient \( u_{k}\) of \( u\) by 
\( u_{k}=\int^{1}_{0}u(\xi ){\rm e}^{-2\pi {\rm i} k\xi}{\rm d}\xi \), where \( {\rm i}\) is the imaginary unit.

For \( 0<T< \infty \) and a Banach space \( X\), 
\( C^{m}([0,T];X) \) denotes the space of functions that are \( m\) times continuously differentiable 
in \( t\) with respect to the norm of \( X\).

For any function space described above, we say that a vector valued function belongs to the function space 
if each of its components does. 
Additionally, for a vector valued function \( \mbox{\mathversion{bold}$u$}\), 
the \( k\)-th Fourier coefficient of \( \mbox{\mathversion{bold}$u$}\) is 
understood as being the vector composed of the \( k\)-th Fourier coefficient of each component of 
\( \mbox{\mathversion{bold}$u$}\).

\vspace*{5mm}

Now we state our main theorem regarding the solvability of (\ref{ext}). 
\begin{Th}\label{TH}
Let \( T>0\) and $m$ an integer satisfying \( m\geq 5\). 
If the initial vortex filament \( \mbox{\mathversion{bold}$x$}_{0}\) 
and the external flow \( \mbox{\mathversion{bold}$F$}\) 
satisfy \( \mbox{\mathversion{bold}$x$}_{0} \in H^{m}(\mathbf{T})\), 
\( \inf_{\xi\in\mathbf{T}}|\mbox{\mathversion{bold}$x$}_{0\xi}(\xi)| > 0 \), and 
\( \mbox{\mathversion{bold}$F$}\in C\big( [0,T];W^{m,\infty}(\mathbf{R}^{3})\big)\),
then there exists \( T_{0}\in (0,T]\) such that the initial value problem \eqref{ext} has a unique solution 
\( \mbox{\mathversion{bold}$x$}(\xi ,t) \) in the class 
\( \mbox{\mathversion{bold}$x$} \in 
C\big( [0,T_{0}]; H^{m}(\mathbf{T})\big)\cap C^{1}\big( [0,T_{0}]; H^{m-2}(\mathbf{T})\big) \) and 
\( \inf_{\xi\in\mathbf{T}}|\mbox{\mathversion{bold}$x$}_{\xi}(\xi,t)| > 0 \) for \( 0\leq t\leq T_{0} \).
\end{Th}

The above theorem gives the time-local solvability of (\ref{ext}). We note that 
Nishiyama \cite{10} proved the existence of the solution in \( C\big( [0,T]; H^{2}(\mathbf{T})\big) \)
for any \( T>0\), but the uniqueness was not shown. Our result is an extension of his result in that 
we prove the solvability in a more regular Sobolev space together with the uniqueness of the solution. 
The rest of the paper is devoted to the proof of Theorem \ref{TH}.

%
%
%
%
%
\section{Construction of the Solution}
\setcounter{equation}{0}
In this section, we give a brief explanation regarding the construction of the solution. 
The method shown in this section is due to Nishiyama \cite{10}. 
We construct the solution to problem (\ref{ext}) by passing to the limit \( \varepsilon \rightarrow +0\)
in the following regularized problem.
\begin{eqnarray}
\left\{
\begin{array}{ll}
\displaystyle \mbox{\mathversion{bold}$x$}_{t} 
= -\varepsilon \mbox{\mathversion{bold}$x$}_{\xi \xi \xi \xi}
+\frac{\mbox{\mathversion{bold}$x$}_{\xi}\times \mbox{\mathversion{bold}$x$}_{\xi\xi}}
{|\mbox{\mathversion{bold}$x$}_{\xi}|^{3}+\varepsilon ^{\alpha}} 
+ \mbox{\mathversion{bold}$F$}( \mbox{\mathversion{bold}$x$},t),& \xi \in \mathbf{T}, \ t>0, \\[2ex]
\mbox{\mathversion{bold}$x$}(\xi, 0) = \mbox{\mathversion{bold}$x$}_{0}(\xi), & \xi \in \mathbf{T},
\end{array}\right.
\label{regnl}
\end{eqnarray}
where \( \varepsilon >0\) and \( \alpha \) with \( 0<\alpha <3/8 \) are real parameters. 
The solution of problem (\ref{regnl}) can be constructed by an iteration scheme based on the 
solvability of the following linear problem.
\begin{eqnarray}
\left\{
\begin{array}{ll}
\displaystyle \mbox{\mathversion{bold}$x$}_{t} = 
-\varepsilon \mbox{\mathversion{bold}$x$}_{\xi \xi \xi \xi } + \mbox{\mathversion{bold}$G$}, 
& \xi \in \mathbf{T}, \ t>0, \\[1ex]
\mbox{\mathversion{bold}$x$}(\xi ,0) = \mbox{\mathversion{bold}$x$}_{0}(\xi ), & \xi \in \mathbf{T}.
\end{array}\right.
\label{linear}
\end{eqnarray}
Finally, 
for \( \mbox{\mathversion{bold}$x$}_{0}\in H^{m}(\mathbf{T})\) and 
\( \mbox{\mathversion{bold}$G$}\in C\big([0,T];W^{m-2,\infty}(\mathbf{T})\big) \),
the unique existence of the solution to (\ref{linear}) in 
\( C\big( [0,T]; H^{m}(\mathbf{T})\big)\cap C^{1}\big( [0,T]; H^{m-4}(\mathbf{T})\big) \)
for any \( T>0\) and \( m\geq 4\) is known from the standard theory of parabolic equations.
In fact, multiplying the first equation in (\ref{linear}) by \( {\rm e}^{2\pi {\rm i} k \xi }\) 
for \( k\in \mathbf{Z}\) 
and integrating with respect to \( \xi \) over \( \mathbf{T}\), we see that
the solution \( \mbox{\mathversion{bold}$x$}(\xi ,t) \) of (\ref{linear}) is 
given by 
\( \mbox{\mathversion{bold}$x$}(\xi,t)=
\sum_{k\in\mathbf{Z}}\mbox{\mathversion{bold}$x$}_{k}(t){\rm e}^{2\pi {\rm i}k\xi} \), 
where \( \mbox{\mathversion{bold}$x$}_{k}\) is the solution of the following ordinary differential equation. 
\begin{align*}
\left\{
\begin{array}{ll}
\displaystyle 
\frac{{\rm d} \mbox{\mathversion{bold}$x$}_{k}}{{\rm d} t}
=
-16\pi^{4}k^{4}\varepsilon\mbox{\mathversion{bold}$x$}_{k} + \mbox{\mathversion{bold}$G$}_{k},
& t>0, \\[1ex]
\mbox{\mathversion{bold}$x$}_{k}(0)=\mbox{\mathversion{bold}$x$}_{0,k}. & \ 
\end{array}\right.
\end{align*}
Here, \( \mbox{\mathversion{bold}$x$}_{0,k}\) and 
\(\mbox{\mathversion{bold}$G$}_{k} \) are the 
\( k\)-th Fourier coefficients of \( \mbox{\mathversion{bold}$x$}_{0}\) and 
\( \mbox{\mathversion{bold}$G$}\), respectively.
The solution of the above ordinary differential equation is 
given explicitly by
\begin{align*}
\mbox{\mathversion{bold}$x$}_{k}(t)
=
{\rm e}^{-16\pi^{4}k^{4}\varepsilon t}\mbox{\mathversion{bold}$x$}_{0,k}
+
\int^{t}_{0}
{\rm e}^{-16\pi^{4}k^{4}\varepsilon (t- \tau )}
\mbox{\mathversion{bold}$G$}_{k}(\tau ) {\rm d}\tau .
\end{align*}
The explicit form of the solution to (\ref{linear}) and direct calculations 
utilizing Parseval's equality show that
the sequence 
\( \{\mbox{\mathversion{bold}$x$}^{(n)}\}_{n=1}^{\infty}\)
given for \( n \geq 2\) as the solution of 
\begin{align*}
\left\{
\begin{array}{ll}
\displaystyle \mbox{\mathversion{bold}$x$}^{(n)}_{t} 
= -\varepsilon \mbox{\mathversion{bold}$x$}^{(n)}_{\xi \xi \xi \xi}
+\frac{\mbox{\mathversion{bold}$x$}^{(n-1)}_{\xi}
\times \mbox{\mathversion{bold}$x$}^{(n-1)}_{\xi\xi}}
{|\mbox{\mathversion{bold}$x$}^{(n-1)}_{\xi}|^{3}+\varepsilon ^{\alpha}} 
+ \mbox{\mathversion{bold}$F$}( \mbox{\mathversion{bold}$x$}^{(n-1)},t),
& \xi \in \mathbf{T}, \ t>0, \\[2ex]
\mbox{\mathversion{bold}$x$}^{(n)}(\xi, 0) 
= \mbox{\mathversion{bold}$x$}_{0}(\xi), & \xi \in \mathbf{T},
\end{array}\right.
\end{align*}
with \( \mbox{\mathversion{bold}$x$}^{(1)}=\mbox{\mathversion{bold}$x$}_{0}\)
converges to the solution of (\ref{regnl}) in the desired function space.
It is shown in \cite{10} that a solution of (\ref{regnl}) belonging to 
\( C\big( [0,T]; H^{2}(\mathbf{T})\big) \) satisfies 
\( |\mbox{\mathversion{bold}$x$}_{\xi}(\xi ,t)|\geq c_1>0\) for some positive constant
\( c_1\) for all \( \xi \in \mathbf{T}\) and \( t\in [0,T]\). 

We state the existence theorem for convenience.

\begin{Th}
Let \( T>0 \) and \(c_0>0\). There exists \( \varepsilon_{0}>0\) and \( c_1>0\) such that
for any \( \varepsilon \in (0,\varepsilon_{0}]\), integer \( m\geq 4\), and \(0<\alpha <\frac{3}{8}\), 
if the initial vortex filament \( \mbox{\mathversion{bold}$x$}_{0} \) and the external flow 
\( \mbox{\mathversion{bold}$F$} \) satisfy \( \mbox{\mathversion{bold}$x$}_{0}\in H^{m}(\mathbf{T})\), 
\( \inf_{\xi\in\mathbf{T}}|\mbox{\mathversion{bold}$x$}_{0\xi}(\xi)|\geq c_0\), and
\( \mbox{\mathversion{bold}$F$}\in C\big( [0,T];W^{m,\infty}(\mathbf{R}^{3})\big)\), then there exists a 
unique solution \(  \mbox{\mathversion{bold}$x$}(\xi ,t) \) to {\rm (\ref{regnl})} satisfying
\( \mbox{\mathversion{bold}$x$}\in 
C\big( [0,T]; H^{m}(\mathbf{T})\big)\cap C^{1}\big( [0,T]; H^{m-2}(\mathbf{T})\big) \)
and \( |\mbox{\mathversion{bold}$x$}_{\xi}(\xi ,t)|\geq c_1\) for all 
\( \xi \in \mathbf{T}\) and \( t\in [0,T]\).

\end{Th}


\section{Energy Estimates of the Solution}
\setcounter{equation}{0}

Our next step is to derive energy estimates for the solution to (\ref{regnl}) 
which are uniform with respect to \( \varepsilon >0 \).
This will allow us to pass to the limit \( \varepsilon \rightarrow +0\) and finish the proof of 
Theorem \ref{TH}.
We do this by deriving suitable energies 
that allow us to estimate the solution in the function space stated in the theorem. 
The derivation of such energy is the most important part of the proof and thus, we go into more detail.
For simplicity, we derive energy estimates for the solution to our original problem
(\ref{ext}) because the arguments for the uniform estimates of the solution to 
(\ref{regnl}) are the same. 

Thus, our objective is to derive energy estimates for the solution of
\begin{eqnarray}
\left\{
\begin{array}{ll}
\displaystyle \mbox{\mathversion{bold}$x$}_{t} = \frac{\mbox{\mathversion{bold}$x$}_{\xi }\times 
\mbox{\mathversion{bold}$x$}_{\xi \xi}}{|\mbox{\mathversion{bold}$x$}_{\xi}|^{3}}
+ \mbox{\mathversion{bold}$F$}(\mbox{\mathversion{bold}$x$},t), &
\xi \in \mathbf{T}, \ t>0, \\[2ex]
\mbox{\mathversion{bold}$x$}(\xi ,0) = \mbox{\mathversion{bold}$x$}_{0}(\xi), & \xi \in \mathbf{T},
\end{array}\right.
\label{ext2}
\end{eqnarray}
belonging to \( C\big( [0,T]; H^{m}(\mathbf{T})\big)\cap C^{1}\big( [0,T]; H^{m-2}(\mathbf{T})\big) \)
on some time interval \( [0,T_{0}]\) with \( T_{0}\in (0,T] \).
The difficulty arises from the fact that a solution of (\ref{ext2}) stretches, i.e., 
\( |\mbox{\mathversion{bold}$x$}_{\xi}|  \not\equiv 1\) even if 
\( |\mbox{\mathversion{bold}$x$}_{0 \xi}|\equiv 1\). When \( |\mbox{\mathversion{bold}$x$}_{\xi}|\equiv 1\),
many useful properties of the solution can be utilized to obtain energy estimates, but 
these properties are not at our disposal in the present problem setting. 
To overcome this, we modify the energy function from the usual Sobolev norm 
to derive the necessary estimates.

In the following we will derive an a priori estimate under the assumption that the solution \( \boldsymbol{x} \) 
to \eqref{ext2} satisfies 
\begin{equation}\label{hyp1}
 |\boldsymbol{x}_{\xi}(\xi,t)| \geq c_1, \qquad
 \|\boldsymbol{x}(t)\|_{W^{3,\infty}(\mathbf{T})} \leq M_1
\end{equation}
for any \( \xi\in\mathbf{T} \) and any \( t\in[0,T_0] \), where positive constants \( c_1,M_1 \), and 
\(T_0 \in(0,T] \) will be determined later. 
We will denote constants depending on these constants $c_1$ and $M_1$ by the same symbol $C_1$, 
which may change from line to line. 

First, we set \( \mbox{\mathversion{bold}$v$} := 
\mbox{\mathversion{bold}$x$}_{\xi}\) and take the \( \xi \) derivative of 
(\ref{ext2}) to rewrite the equation in terms of \( \mbox{\mathversion{bold}$v$}\).
\begin{eqnarray}
\left\{
\begin{array}{ll}
\mbox{\mathversion{bold}$v$}_{t} = f\mbox{\mathversion{bold}$v$}\times \mbox{\mathversion{bold}$v$}_{\xi \xi}
+ f_{\xi}\mbox{\mathversion{bold}$v$}\times \mbox{\mathversion{bold}$v$}_{\xi}
+ ( {\rm D}\mbox{\mathversion{bold}$F$})\mbox{\mathversion{bold}$v$}, & \xi \in \mathbf{T}, \ t>0, \\[1ex]
\mbox{\mathversion{bold}$v$}(\xi,0) = \mbox{\mathversion{bold}$v$}_{0}(\xi ), & \xi \in \mathbf{T},
\end{array}\right.
\label{v}
\end{eqnarray}
where we have set \( \mbox{\mathversion{bold}$v$}_{0}:= \mbox{\mathversion{bold}$x$}_{0\xi }\),
\( f:=1/|\mbox{\mathversion{bold}$v$}|^{3}\), and omitted the arguments of 
\( \mbox{\mathversion{bold}$F$}\). 
Following standard procedures, we differentiate the equation in \eqref{v} \( k\) times with respect to \( \xi \) 
and set \( \mbox{\mathversion{bold}$v$}^{k} := \partial ^{k}_{\xi}\mbox{\mathversion{bold}$v$}\) to obtain 
\begin{equation}\label{vk}
\mbox{\mathversion{bold}$v$}^{k}_{t}
= f\mbox{\mathversion{bold}$v$} \times \mbox{\mathversion{bold}$v$}^{k}_{\xi \xi}
 + kf\mbox{\mathversion{bold}$v$}_{\xi} \times \mbox{\mathversion{bold}$v$}^{k}_{\xi}
 + (k+1)f_{\xi}\mbox{\mathversion{bold}$v$} \times \mbox{\mathversion{bold}$v$}^{k}_{\xi}
 - 3f^{5/3}(\boldsymbol{v}\cdot\boldsymbol{v}_\xi^k)\boldsymbol{v}\times\boldsymbol{v}_\xi
 + \mbox{\mathversion{bold}$G$}^{k},
\end{equation}
where \( \mbox{\mathversion{bold}$G$}^{k}\) is the collection of terms that contain \( \xi \)-derivatives 
of \( \mbox{\mathversion{bold}$v$}\) up to order \( k\). 
Moreover, by standard calculus inequalities it satisfies the estimate 
\begin{equation}\label{Gk}
\| \mbox{\mathversion{bold}$G$}^{k}\|_l \leq 
C_1\|\mbox{\mathversion{bold}$v$}\|_{k+l}
\end{equation}
for \( k,l=0,1,2,\ldots \) satisfying $k+l+1 \leq m$. 
Now that we have derived (\ref{vk}), the standard method would be to take the \( L^2(\mathbf{T}) \) inner product of 
\( \mbox{\mathversion{bold}$v$}^{k}\) and (\ref{vk}) to estimate the time evolution of 
\( \| \mbox{\mathversion{bold}$v$}^{k}\| \). 
This is not possible for our equation because the terms with derivatives of 
\( \mbox{\mathversion{bold}$v$}^{k}\) cause a loss of regularity. To avoid such loss, we employ a 
series of change of variables to derive a modified energy function from which we can derive the necessary estimates. 
The key idea is to decompose \( \mbox{\mathversion{bold}$v$}^{k}\) into two parts. 
More precisely, we decompose \( \mbox{\mathversion{bold}$v$}^{k}\) as
\begin{eqnarray}
\mbox{\mathversion{bold}$v$}^{k} = \frac{(\mbox{\mathversion{bold}$v$}
\cdot \mbox{\mathversion{bold}$v$}^{k})}{|\mbox{\mathversion{bold}$v$}|^{2}}\mbox{\mathversion{bold}$v$}
- \frac{1}{|\mbox{\mathversion{bold}$v$}|^{2}}\mbox{\mathversion{bold}$v$}
\times (\mbox{\mathversion{bold}$v$}\times \mbox{\mathversion{bold}$v$}^{k}).
\label{dec}
\end{eqnarray}
The above decomposes \( \mbox{\mathversion{bold}$v$}^{k}\) into the sum of 
its \( \mbox{\mathversion{bold}$v$}\) component and the component orthogonal to 
\( \mbox{\mathversion{bold}$v$}\). The decomposition is well-defined since 
we have \( |\mbox{\mathversion{bold}$v$}(\xi,t)|\geq c_1>0\) under our hypotheses. 
The principle part of the components are 
\( \mbox{\mathversion{bold}$v$}\cdot \mbox{\mathversion{bold}$v$}^{k}\) and 
\( \mbox{\mathversion{bold}$v$}\times \mbox{\mathversion{bold}$v$}^{k}\) respectively, and we define
two new variables
\begin{equation}
h^{k}:= \mbox{\mathversion{bold}$v$}\cdot \mbox{\mathversion{bold}$v$}^{k}, \qquad
\mbox{\mathversion{bold}$z$}^{k}:= \mbox{\mathversion{bold}$v$}\times \mbox{\mathversion{bold}$v$}^{k},
\end{equation}
and estimate them separately.


\subsection{Estimate of \( h^{k}\)}
\label{hk}
We first derive an equation for \( h^{k}\). 
If follows from \eqref{vk} that 
\begin{align*}
h_t^k
&= \boldsymbol{v}\cdot\boldsymbol{v}_t^k + \boldsymbol{v}_t\cdot\boldsymbol{v}^k \\
&= kf(\boldsymbol{v}\times\boldsymbol{v}_{\xi})\cdot\boldsymbol{v}_{\xi}^k
 + \boldsymbol{v}\cdot\boldsymbol{G}^k + \boldsymbol{v}_t\cdot\boldsymbol{v}^k
\end{align*}
and that 
\begin{align*}
(\boldsymbol{v}_{\xi}\cdot\boldsymbol{v}^k)_t
&= \boldsymbol{v}_{\xi}\cdot\boldsymbol{v}_t^k + \boldsymbol{v}_{t\xi}\cdot\boldsymbol{v}^k \\
&= -f(\boldsymbol{v}\times\boldsymbol{v}_{\xi})\cdot\boldsymbol{v}_{\xi\xi}^k
 + \boldsymbol{v}_{\xi}\cdot( (k+1)f_\xi\boldsymbol{v}\times\boldsymbol{v}_\xi^k + \boldsymbol{G}^k )
  + \boldsymbol{v}_{t\xi}\cdot\boldsymbol{v}^k,
\end{align*}
which yield 
\begin{align*}
( h_\xi^k + k\boldsymbol{v}_{\xi}\cdot\boldsymbol{v}^k )_t 
&= k(f(\boldsymbol{v}\times\boldsymbol{v}_{\xi}))_{\xi}\cdot\boldsymbol{v}_{\xi}^k
 + ( \boldsymbol{v}\cdot\boldsymbol{G}^k + \boldsymbol{v}_t\cdot\boldsymbol{v}^k )_{\xi} \\
&\quad\;
 + k\{ \boldsymbol{v}_{\xi}\cdot( (k+1)f_\xi\boldsymbol{v}\times\boldsymbol{v}_\xi^k + \boldsymbol{G}^k )
  + \boldsymbol{v}_{t\xi}\cdot\boldsymbol{v}^k \} \\
&=: \boldsymbol{G}_1^k.
\end{align*}
Here, if we impose an additional assumption $2\leq k\leq m-2$, then by \eqref{Gk} and the standard calculus 
inequalities we have 
\[
\|\boldsymbol{G}_1^k\| \leq C_1\|\boldsymbol{v}\|_{k+1},
\]
where we used the equation in \eqref{v} to replace the $t$-derivative with $\xi$-derivatives. 
Therefore, we obtain 
\begin{equation}\label{estimate hk}
\frac{{\rm d}}{{\rm d}t} \|h_\xi^k + k\boldsymbol{v}_{\xi}\cdot\boldsymbol{v}^k\|^2
 \leq C_1\|\boldsymbol{v}\|_{k+1}^2
\end{equation}
for \( 2\leq k\leq m-2 \).

%

\subsection{Estimate of \( \mbox{\mathversion{bold}$z$}^{k}\)}
\label{Zk}
Next we consider \( \mbox{\mathversion{bold}$z$}^{k}\) and derive an equation for \( \boldsymbol{z}^k \). 
It follows from \eqref{vk} and the identity 
\( \boldsymbol{z}_{\xi\xi}^k = \boldsymbol{v}\times\boldsymbol{v}_{\xi\xi}^k
 + 2\boldsymbol{v}_\xi\times\boldsymbol{v}_{\xi}^k + \boldsymbol{v}_{\xi\xi}\times\boldsymbol{v}^k$ that 
\begin{align*}
\boldsymbol{z}_t^k
&= \boldsymbol{v}\times\{ f( \boldsymbol{z}_{\xi\xi}^k + (k-2)\boldsymbol{v}_\xi\times\boldsymbol{v}_{\xi}^k
 - \boldsymbol{v}_{\xi\xi}\times\boldsymbol{v}^k ) \\
&\qquad
 + (k+1)f_\xi\boldsymbol{v}\times\boldsymbol{v}_\xi^k
 - 3f^{5/3}(\boldsymbol{v}\cdot\boldsymbol{v}_\xi^k)\boldsymbol{v}\times\boldsymbol{v}_\xi + \boldsymbol{G}^k \}
 + \boldsymbol{v}_t\times\boldsymbol{v}^k.
\end{align*}
Here, by using the decomposition \eqref{dec} with $\boldsymbol{v}^k$ replaced by $\boldsymbol{v}_\xi^k$ 
and $\boldsymbol{v}_\xi$ we have 
\begin{align}\label{rel1}
\boldsymbol{v}\times(\boldsymbol{v}_\xi\times\boldsymbol{v}_\xi^k)
&= (\boldsymbol{v}\cdot\boldsymbol{v}_\xi^k)\boldsymbol{v}_\xi
 - (\boldsymbol{v}\cdot\boldsymbol{v}_\xi)\boldsymbol{v}_\xi^k \\
&= (\boldsymbol{v}\cdot\boldsymbol{v}_\xi^k)\boldsymbol{v}_\xi
 - (\boldsymbol{v}\cdot\boldsymbol{v}_\xi)\biggl(
  \frac{(\boldsymbol{v}\cdot\boldsymbol{v}_\xi^k)}{|\boldsymbol{v}|^2}\boldsymbol{v}
   - \frac{1}{|\boldsymbol{v}|^2}\boldsymbol{v}\times(\boldsymbol{v}\times\boldsymbol{v}_\xi^k) \biggr) \nonumber \\
&= \biggl( \boldsymbol{v}_\xi - \frac{(\boldsymbol{v}\cdot\boldsymbol{v}_\xi)}{|\boldsymbol{v}|^2}\boldsymbol{v}
    \biggr)(\boldsymbol{v}\cdot\boldsymbol{v}_\xi^k)
 + \frac{(\boldsymbol{v}\cdot\boldsymbol{v}_\xi)}{|\boldsymbol{v}|^2}
  \boldsymbol{v}\times(\boldsymbol{v}\times\boldsymbol{v}_\xi^k) \nonumber \\
&= -\frac{1}{|\boldsymbol{v}|^2}
 (\boldsymbol{v}\times(\boldsymbol{v}\times\boldsymbol{v}_\xi))(\boldsymbol{v}\cdot\boldsymbol{v}_\xi^k)
 - \frac13f^{-1}f_\xi\boldsymbol{v}\times(\boldsymbol{v}\times\boldsymbol{v}_\xi^k). \nonumber 
\end{align}
These together with the identity \( \boldsymbol{z}_{\xi}^k = \boldsymbol{v}\times\boldsymbol{v}_{\xi}^k
 + \boldsymbol{v}_{\xi}\times\boldsymbol{v}^k$ yield 
\begin{equation}\label{eqzk}
\boldsymbol{z}_t^k
= f\boldsymbol{v}\times\bigl( \boldsymbol{z}_{\xi\xi}^k
 - (k+1)f^{2/3}(\boldsymbol{v}\times\boldsymbol{v}_\xi)(\boldsymbol{v}\cdot\boldsymbol{v}_\xi^k) \bigr)
 + \frac{2k+5}{3}f_\xi\boldsymbol{v}\times\boldsymbol{z}_\xi^k + \boldsymbol{G}_2^k,
\end{equation}
where 
\[
\boldsymbol{G}_2^k
= -\frac{2k+5}{3}f_\xi\boldsymbol{v}_\xi\times\boldsymbol{v}^k
 + \boldsymbol{v}\times( \boldsymbol{G}^k-\boldsymbol{v}_{\xi\xi}\times\boldsymbol{v}^k )
 + \boldsymbol{v}_t\times\boldsymbol{v}^k,
\]
which satisfies the estimate 
\[
\| \boldsymbol{G}_2^k \|_1 \leq C_1\|\boldsymbol{v}\|_{k+1}
\]
for $2\leq k\leq m-2$. 
Here, we used \eqref{Gk} and the standard calculus inequalities together with the equation in \eqref{v} 
to replace the $t$-derivative with $\xi$-derivatives.

In view of \eqref{eqzk}, we make a change of variable given by 
\begin{equation}\label{defuk}
\boldsymbol{u}^k := \boldsymbol{z}^k
 - (k+1)f^{2/3}(\boldsymbol{v}\times\boldsymbol{v}_\xi)( \boldsymbol{v}\cdot\boldsymbol{v}^{k-1} ). 
\end{equation}
It follows from \eqref{vk} with \( k \) replacaed by \( k-1 \) that 
\( \boldsymbol{v}\cdot\boldsymbol{v}_t^{k-1}
 = \boldsymbol{v}\cdot( (k-1)f\boldsymbol{v}_\xi\times\boldsymbol{v}^k + \boldsymbol{G}^{k-1} ) \), 
which together with \eqref{eqzk} implies 
\begin{align*}
\boldsymbol{u}_t^k
&= f\boldsymbol{v}\times\bigl( \boldsymbol{z}_{\xi\xi}^k
 - (k+1)f^{2/3}(\boldsymbol{v}\times\boldsymbol{v}_\xi)(\boldsymbol{v}\cdot\boldsymbol{v}_\xi^k) \bigr)
 + \frac{2k+5}{3}f_\xi\boldsymbol{v}\times\boldsymbol{z}_\xi^k + \boldsymbol{G}_2^k \\
&\quad\;
 - (k+1)f^{2/3}(\boldsymbol{v}\times\boldsymbol{v}_\xi)\bigl( \boldsymbol{v}\cdot
  ( (k-1)f\boldsymbol{v}_\xi\times\boldsymbol{v}^k + \boldsymbol{G}^{k-1} ) \bigr) \\
&\quad\;
 - (k+1)( f^{2/3}(\boldsymbol{v}\times\boldsymbol{v}_\xi)\otimes\boldsymbol{v} )_t\boldsymbol{v}^{k-1}.
\end{align*}
Differentiating \eqref{defuk} with respect to $\xi$ we have 
\begin{align*}
\boldsymbol{z}_\xi^k 
&= \boldsymbol{u}_\xi^k
 + (k+1)(f^{2/3}(\boldsymbol{v}\times\boldsymbol{v}_\xi)( \boldsymbol{v}\cdot\boldsymbol{v}^{k-1} ))_\xi, \\
\boldsymbol{z}_{\xi\xi}^k 
&= \boldsymbol{u}_{\xi\xi}^k
 + (k+1)f^{2/3}(\boldsymbol{v}\times\boldsymbol{v}_\xi)( \boldsymbol{v}\cdot\boldsymbol{v}_\xi^k ) \\
&\quad\;
 + (k+1)\{ 2(f^{2/3}(\boldsymbol{v}\times\boldsymbol{v}_\xi)\otimes\boldsymbol{v})_{\xi}\boldsymbol{v}^k
  + (f^{2/3}(\boldsymbol{v}\times\boldsymbol{v}_\xi)\otimes\boldsymbol{v})_{\xi\xi}\boldsymbol{v}^{k-1} \}.
\end{align*}
Substituting these into the above equation, we obtain 
\begin{equation}\label{equk}
\boldsymbol{u}_t^k = f\boldsymbol{v}\times\boldsymbol{u}_{\xi\xi}^k
 + \frac{2k+5}{3}f_\xi\boldsymbol{v}\times\boldsymbol{u}_\xi^k + \boldsymbol{G}_3^k,
\end{equation}
where 
\begin{align*}
\boldsymbol{G}_3^k
&= (k+1)f\boldsymbol{v}\times\{ 
 2(f^{2/3}(\boldsymbol{v}\times\boldsymbol{v}_\xi)\otimes\boldsymbol{v})_{\xi}\boldsymbol{v}^k
  + (f^{2/3}(\boldsymbol{v}\times\boldsymbol{v}_\xi)\otimes\boldsymbol{v})_{\xi\xi}\boldsymbol{v}^{k-1} \} \\
&\quad\;
 + \frac{2k+5}{3}(k+1)f_\xi\boldsymbol{v}\times
  (f^{2/3}(\boldsymbol{v}\times\boldsymbol{v}_\xi)( \boldsymbol{v}\cdot\boldsymbol{v}^{k-1} ))_\xi
 + \boldsymbol{G}_2^k \\
&\quad\;
 - (k+1)f^{2/3}(\boldsymbol{v}\times\boldsymbol{v}_\xi)\bigl( \boldsymbol{v}\cdot
  ( (k-1)f\boldsymbol{v}_\xi\times\boldsymbol{v}^k + \boldsymbol{G}^{k-1} ) \bigr) \\
&\quad\;
 - (k+1)( f^{2/3}(\boldsymbol{v}\times\boldsymbol{v}_\xi)\otimes\boldsymbol{v} )_t\boldsymbol{v}^{k-1}.
\end{align*}
If we further impose an additional assumption \( 3\leq k\leq m-2 \), then as before we have
\[
\| \boldsymbol{G}_3^k \|_1 \leq C_1\|\boldsymbol{v}\|_{k+1}. 
\]

We further make a change of variable, which could be perceived as a type of gauge transformation,
to negate the loss of regularity caused by the second term containing 
\( \mbox{\mathversion{bold}$u$}^k_{\xi}\) on the right-hand side of \eqref{equk}.
We do this by changing the variable from \( \mbox{\mathversion{bold}$u$}^k\) to 
\( \mbox{\mathversion{bold}$w$}^k\) in the form \( \mbox{\mathversion{bold}$u$}^k=
a(\xi,t)\mbox{\mathversion{bold}$w$}^k\) for some positive scalar function \( a(\xi,t)\),
that is harmless to our energy estimate, to cancel out the terms causing the loss of regularity. 
Substituting this into \eqref{equk} yields 
\begin{align*}
a\boldsymbol{w}^k_t + a_{t}\boldsymbol{w}^k
&= af\boldsymbol{v}\times \boldsymbol{w}^k_{\xi\xi} 
 + \biggl( 2fa_\xi + \frac{2k+5}{3}f_\xi a \biggr) \boldsymbol{v}\times\boldsymbol{w}_\xi^k \\
&\quad\;
 + \biggl( fa_{\xi\xi} + \frac{2k+5}{3}f_\xi a_\xi \biggr) \boldsymbol{v}\times\boldsymbol{w}^k
 + \boldsymbol{G}_3^k.
\end{align*}
Hence, if we can choose \( a(\xi,t)\) so that
\[
2fa_\xi + \frac{2k+5}{3}f_\xi a=0,
\]
then the terms causing the loss of regularity are canceled. 
Dividing the above relation by \( 2af\) yields 
\( \bigl( \log(af^{\frac{2k+5}{6}}) \bigr)_\xi =0 \). 
Therefore, it is sufficient to choose \( a(\xi,t)\) by 
\begin{equation}\label{a}
a(\xi,t) := f(\xi,t)^{-\frac{2k+5}{6}} = |\boldsymbol{v}(\xi,t)|^{k+\frac{5}{2}}. 
\end{equation}
Then, the equation for \( \boldsymbol{w}^k \) becomes 
\begin{equation}\label{eqwk}
\boldsymbol{w}^k_t = f\boldsymbol{v}\times \boldsymbol{w}^k_{\xi\xi} + \boldsymbol{G}_4^k,
\end{equation}
where 
\[
\boldsymbol{G}_4^k = a^{-1}\biggl\{ 
 \biggl( fa_{\xi\xi} + \frac{2k+5}{3}f_\xi a_\xi \biggr) \boldsymbol{v}\times\boldsymbol{w}^k
 + \boldsymbol{G}_3^k - a_{t}\boldsymbol{w}^k \biggr\}.
\]
As before, under the condition \( 3\leq k\leq m-2 \) we have 
\[
\| \boldsymbol{G}_4^k \|_1 \leq C_1\| \boldsymbol{v} \|_{k+1}.
\]
Now, taking the \( \xi \) derivative of \eqref{eqwk} yields 
\[
\boldsymbol{w}^k_{\xi t} = (f\boldsymbol{v}\times \boldsymbol{w}^k_{\xi\xi})_\xi + \boldsymbol{G}_{4\xi}^k,
\]
from which we can easily derive the estimate 
\begin{equation}\label{wk}
\frac{{\rm d}}{{\rm d}t} \| \boldsymbol{w}_\xi^k \|^2 \leq C_1\|\boldsymbol{v}\|_{k+1}^2
\end{equation}
for \( 3\leq k\leq m-2 \).

%

\subsection{Completion of a priori estimate}
We suppose \( m\geq5 \) and derive an a priori estimate of the solution $\boldsymbol{x}$ 
to the initial value problem \eqref{ext2} under the hypotheses 
\begin{equation}\label{hyp2}
|\boldsymbol{x}_{\xi}(\xi,t)| \geq c_1, \quad
 \|\boldsymbol{x}(t)\|_{W^{2,\infty}(\mathbf{T})} \leq M_0, \quad
 \|\boldsymbol{x}(t)\|_m \leq M_1
\end{equation}
for any \( \xi\in\mathbf{T} \) and any \( t\in[0,T_0] \), where positive constants \( c_1, M_0, M_1 \), 
and \(T_0 \in (0,T] \) will be determined later. 
These hypotheses will be justified later as usual. 
We also note that the last inequality in \eqref{hyp2} and the Sobolev imbedding theorem give 
\( \|\boldsymbol{x}(t)\|_{W^{3,\infty}(\mathbf{T})} \leq CM_1 \) with an absolute constant $C$. 
In the following we simply write the constants depending on $c_1$ and $M_0$ by $C_0$ and the 
constants depending also on $M_1$ by $C_1$ as before, which may change from line to line.

It follows directly from \eqref{ext2} that 
\begin{equation}\label{esx}
\frac{\rm d}{{\rm d}t}\|\boldsymbol{x}(t)\|^2 \leq C_0.
\end{equation}
Now, we define a modified energy function \( E^k(t) \) by 
\[
E^k(t) := \| h_\xi^k(t) + \boldsymbol{v}_\xi(t)\cdot\boldsymbol{v}^k(t) \|^2
 + \| \boldsymbol{w}_\xi^k(t) \|^2 + \| \boldsymbol{x}(t) \|^2.
\]
By the decomposition \eqref{dec} with $\boldsymbol{v}^k$ replaced by $\boldsymbol{v}^{k+1}$ and 
the identities \( h_\xi^k=\boldsymbol{v}\cdot\boldsymbol{v}^{k+1}+\boldsymbol{v}_\xi\cdot\boldsymbol{v}^k \) and 
\( \boldsymbol{z}_\xi^k=\boldsymbol{v}\times\boldsymbol{v}^{k+1}+\boldsymbol{v}_\xi\times\boldsymbol{v}^k \), 
we have 
\begin{align*}
\boldsymbol{v}^{k+1}
&= f^{2/3}(\boldsymbol{v}\cdot\boldsymbol{v}^{k+1})\boldsymbol{v}
 - f^{2/3}\boldsymbol{v}\times(\boldsymbol{v}\times\boldsymbol{v}^{k+1}) \\
&= f^{2/3}(h_\xi^k+k\boldsymbol{v}_\xi\cdot\boldsymbol{v}^k)\boldsymbol{v}
 - af^{2/3}\boldsymbol{v}\times\boldsymbol{w}^k_\xi + \boldsymbol{G}_5^k,
\end{align*}
where 
\[
\boldsymbol{G}_5^k
= -(k+1)f^{2/3}\boldsymbol{v}_\xi\cdot\boldsymbol{v}^k
 -f^{2/3}\boldsymbol{v}\times\{a_\xi\boldsymbol{w}^k
  + (k+1)( f^{2/3}(\boldsymbol{v}\times\boldsymbol{v}_\xi)(\boldsymbol{v}\cdot\boldsymbol{v}^{k-1}) )_\xi
  -\boldsymbol{v}_\xi\times\boldsymbol{v}^k \},
\]
which satisfies 
\( \| \boldsymbol{G}_5^k \| \leq C_0 \|\boldsymbol{v}^k\| \). 
Therefore, in view of \( \|\boldsymbol{x}\|_{k+2} \leq C_0( \|\boldsymbol{x}\| + \|\boldsymbol{v}^{k+1}\| ) \) 
and the interpolation inequality 
\( \|\boldsymbol{v}^k\| \leq \epsilon\|\boldsymbol{v}^{k+1}\| + C_\epsilon\|\boldsymbol{x}\| \) for \(\epsilon>0\), 
we obtain the equivalence 
\begin{equation}\label{equi}
C_0^{-1} \|\boldsymbol{x}(t)\|_{k+2}^2 \leq E^k(t) \leq C_0 \|\boldsymbol{x}(t)\|_{k+2}^2. 
\end{equation}
Now, we choose \( k=m-2 \). 
Adding \eqref{estimate hk}, \eqref{wk}, and \eqref{esx}, we have 
\[
\frac{\rm d}{{\rm d}t}E^k(t) \leq C_1\|\boldsymbol{v}\|_{k+1}^2 + C_0 \leq C_1( E^k(t) + 1 ),
\]
so that Gronwall's inequality and \eqref{equi} yield 
\[
\|\boldsymbol{x}(t)\|_m^2 \leq C_0\mbox{\rm e}^{C_1t}( \|\boldsymbol{x}_0\|_m^2 + C_1t)
\]
for \( 0\leq t\leq T_0 \).
It follows directly from \eqref{ext2} that 
\( \|\boldsymbol{x}_t(t)\|_{W^{2,\infty}(\mathbf{T})} \leq C_1\) so that we have 
\[
\begin{cases}
 |\boldsymbol{x}_\xi(\xi,t)| \geq |\boldsymbol{x}_{0\xi}(\xi)| - C_1T_0, \\
 \|\boldsymbol{x}(t)\|_{W^{2,\infty}(\mathbf{T})} \leq \|\boldsymbol{x}_0\|_{W^{2,\infty}(\mathbf{T})} + C_1T_0
\end{cases}
\]
for \( 0\leq t\leq T_0 \). 
In view of these estimates, we choose the positive constants $c_1, M_0$, and $M_1$ in \eqref{hyp2} such that 
\( \inf_{\xi\in\mathbf{T}}|\boldsymbol{x}_{0\xi}(\xi)| \geq 2c_1 \), 
\( 2\|\boldsymbol{x}_0\|_{W^{2,\infty}(\boldsymbol{T})} \leq M_0 \), and 
\( 2\sqrt{C_0}\|\boldsymbol{x}_0\|_m \leq M_1\), and then choose $T_0 \in (0,T]$ sufficiently small. 
Then, we see that \eqref{hyp2} holds.


\section{Uniqueness of the Solution}
\setcounter{equation}{0}

The uniqueness of the solution can be proved by the standard method of estimating the difference of two solutions 
with the same initial datum along the line of the energy estimate carried out in the previous section with 
slight modifications. 
In this section, for completeness, we will give the proof.

Suppose \( \mbox{\mathversion{bold}$x$}^{(1)}\) and 
\( \mbox{\mathversion{bold}$x$}^{(2)}\) are solutions of (\ref{ext2}) with the same 
initial datum \( \mbox{\mathversion{bold}$x$}_{0}\) given in Theorem \ref{TH}. 
We set \( \boldsymbol{v}^{(j)} := \mbox{\mathversion{bold}$x$}_\xi^{(j)} \) and 
\( f^{(j)} := 1/|\mbox{\mathversion{bold}$x$}_\xi^{(j)}|^3 \) for \( j=1,2 \) 
and differentiate the equation twice with respect to $\xi$ to obtain 
\begin{align}\label{eqvj}
\boldsymbol{v}_{\xi t}^{(j)}
&= f^{(j)}\boldsymbol{v}^{(j)}\times\boldsymbol{v}_{\xi\xi\xi}^{(j)}
 + f^{(j)}\boldsymbol{v}_\xi^{(j)}\times\boldsymbol{v}_{\xi\xi}^{(j)}
 + 2f_\xi^{(j)}\boldsymbol{v}^{(j)}\times\boldsymbol{v}_{\xi\xi}^{(j)} \\
&\quad\;
 - 3(f^{(j)})^{5/3}(\boldsymbol{v}^{(j)}\cdot\boldsymbol{v}_{\xi\xi}^{(j)})
  (\boldsymbol{v}^{(j)}\times\boldsymbol{v}_{\xi}^{(j)})
 + \boldsymbol{G}^{(j)}, \nonumber
\end{align}
where 
\begin{align*}
\boldsymbol{G}^{(j)}
&= {\rm D}\boldsymbol{F}(\boldsymbol{x}^{(j)},t)\boldsymbol{v}_\xi^{(j)}
 + {\rm D}^2\boldsymbol{F}(\boldsymbol{x}^{(j)},t)[\boldsymbol{v}^{(j)},\boldsymbol{v}^{(j)}] 
 - 3\bigl( ((f^{(j)})^{5/3}\boldsymbol{v}^{(j)})_\xi\cdot\boldsymbol{v}_\xi^{(j)} \bigr)
  (\boldsymbol{v}^{(j)}\times\boldsymbol{v}_{\xi}^{(j)}).
\end{align*}
We also set 
\( \dot{\boldsymbol{x}} := \boldsymbol{x}^{(1)} - \boldsymbol{x}^{(2)} \) and 
\( \dot{\boldsymbol{v}} := \boldsymbol{v}^{(1)} - \boldsymbol{v}^{(2)} = \dot{\boldsymbol{x}}_\xi \), 
and obtain 
\begin{align}\label{eqdv}
\dot{\boldsymbol{v}}_{\xi t}
&= f^{(1)}\boldsymbol{v}^{(1)}\times\dot{\boldsymbol{v}}_{\xi\xi\xi}
 + f^{(1)}\boldsymbol{v}_\xi^{(1)}\times\dot{\boldsymbol{v}}_{\xi\xi}
 + 2f_\xi^{(1)}\boldsymbol{v}^{(1)}\times\dot{\boldsymbol{v}}_{\xi\xi} \\
&\quad\;
 - 3(f^{(1)})^{5/3}(\boldsymbol{v}^{(1)}\cdot\dot{\boldsymbol{v}}_{\xi\xi})
  (\boldsymbol{v}^{(1)}\times\boldsymbol{v}_{\xi}^{(1)})
 + \dot{\boldsymbol{G}}, \nonumber
\end{align}
where 
\begin{align*}
\dot{\boldsymbol{G}}
&= \boldsymbol{G}^{(1)} - \boldsymbol{G}^{(2)}
 + (f^{(1)}\boldsymbol{v}^{(1)} - f^{(2)}\boldsymbol{v}^{(2)})\times\boldsymbol{v}_{\xi\xi\xi}^{(2)} \\
&\quad\;
 + (f^{(1)}\boldsymbol{v}_\xi^{(1)} - f^{(2)}\boldsymbol{v}_\xi^{(2)})\times\boldsymbol{v}_{\xi\xi}^{(2)} 
 + 2(f_\xi^{(1)}\boldsymbol{v}^{(1)} - f_\xi^{(2)}\boldsymbol{v}^{(2)})\times\boldsymbol{v}_{\xi\xi}^{(2)} \\
&\quad\;
 - 3\bigl( (f^{(1)})^{5/3}(\boldsymbol{v}^{(1)}\times\boldsymbol{v}_{\xi}^{(1)})\otimes\boldsymbol{v}^{(1)}
  - (f^{(2)})^{5/3}(\boldsymbol{v}^{(2)}\times\boldsymbol{v}_{\xi}^{(2)})\otimes\boldsymbol{v}^{(2)} \bigr)
  \boldsymbol{v}_{\xi\xi}^{(2)}
\end{align*}
which satisfies 
\( \| \dot{\boldsymbol{G}} \|_1
 \leq C( \|\dot{\boldsymbol{v}}\|_2 + \|\dot{\boldsymbol{x}}\| )
 \leq C\|\dot{\boldsymbol{x}}\|_3\). 
As before, we decompose \( \dot{\boldsymbol{v}}_\xi \) into its \( \boldsymbol{v}^{(1)} \) component 
and the component orthogonal to \( \boldsymbol{v}^{(1)} \) so that we put 
\begin{equation}
\dot{h} := \boldsymbol{v}^{(1)}\cdot\dot{\boldsymbol{v}}_\xi, \qquad
\dot{\boldsymbol{z}} := \boldsymbol{v}^{(1)}\times\dot{\boldsymbol{v}}_\xi
\end{equation}
to obtain 
\( \dot{\boldsymbol{v}}_\xi = (f^{(1)})^{2/3} \bigl( \dot{h}\boldsymbol{v}^{(1)}
 - \boldsymbol{v}^{(1)}\times\dot{\boldsymbol{z}} \bigr) \). 
Then, we have 
\begin{align}
(\dot{h}_\xi+\boldsymbol{v}_\xi^{(1)}\cdot\dot{\boldsymbol{v}}_\xi)_t
&= (f^{(1)}\boldsymbol{v}^{(1)}\times\boldsymbol{v}_\xi^{(1)})_\xi\cdot\dot{\boldsymbol{v}}_{\xi\xi}
 + (\boldsymbol{v}^{(1)}\cdot\dot{\boldsymbol{G}} + \boldsymbol{v}_t^{(1)}\cdot\dot{\boldsymbol{v}}_\xi)_\xi \\
&\quad\;
 + \boldsymbol{v}_\xi^{(1)}\cdot\{
  ( f^{(1)}\boldsymbol{v}_\xi^{(1)} + 2f_\xi^{(1)}\boldsymbol{v}^{(1)} )\times\dot{\boldsymbol{v}}_{\xi\xi} \nonumber \\
&\qquad
  - 3(f^{(1)})^{5/3}(\boldsymbol{v}^{(1)}\times\boldsymbol{v}_\xi^{(1)})
   (\boldsymbol{v}^{(1)}\cdot\dot{\boldsymbol{v}}_{\xi\xi}) + \dot{\boldsymbol{G}} \} \nonumber \\
&=: \dot{\boldsymbol{G}}_1. \nonumber
\end{align}
Here, \( \dot{\boldsymbol{G}}_1 \) satisfies the estimate 
\( \|\dot{\boldsymbol{G}}_1\| \leq C( \|\dot{\boldsymbol{v}}\|_2 + \|\dot{\boldsymbol{x}}\| )
 \leq C\|\dot{\boldsymbol{x}}\|_3 \), 
so that we have 
\begin{equation}\label{estdh}
\frac{\rm d}{{\rm d}t} \| \dot{h}_\xi+\boldsymbol{v}_\xi^{(1)}\cdot\dot{\boldsymbol{v}}_\xi \|^2
 \leq C \|\dot{\boldsymbol{x}}\|_3^2.
\end{equation}

As for \( \dot{\boldsymbol{z}} \) we have 
\begin{align*}
\dot{\boldsymbol{z}}_t
&= f^{(1)}\boldsymbol{v}^{(1)}\times\dot{\boldsymbol{z}}_{\xi\xi}
 - f^{(1)}\boldsymbol{v}^{(1)}\times(\boldsymbol{v}_\xi^{(1)}\times\dot{\boldsymbol{v}}_{\xi\xi}) 
 + 2f_\xi^{(1)}\boldsymbol{v}^{(1)}\times\dot{\boldsymbol{z}}_\xi \\
&\quad\;
 - 3(f^{(1)})^{5/3}(\boldsymbol{v}^{(1)}\times(\boldsymbol{v}^{(1)}\times\boldsymbol{v}_\xi^{(1)}) )
  (\boldsymbol{v}^{(1)}\cdot\dot{\boldsymbol{v}}_{\xi\xi}) \\
&\quad\;
 + \boldsymbol{v}^{(1)}\times( \dot{\boldsymbol{G}}
 - f^{(1)}\boldsymbol{v}_{\xi\xi}^{(1)}\times\dot{\boldsymbol{v}}_\xi
 - 2f_\xi^{(1)}\boldsymbol{v}_\xi^{(1)}\times\dot{\boldsymbol{v}}_\xi )
 + \boldsymbol{v}_t^{(1)}\times\dot{\boldsymbol{v}}_\xi.
\end{align*}
Here, in the same way as the calculations in \eqref{rel1} we have 
\begin{align*}
\boldsymbol{v}^{(1)}\times(\boldsymbol{v}_\xi^{(1)}\times\dot{\boldsymbol{v}}_{\xi\xi})
&= - (f^{(1)})^{2/3}(\boldsymbol{v}^{(1)}\times(\boldsymbol{v}^{(1)}\times\boldsymbol{v}_\xi^{(1)})
  (\boldsymbol{v}^{(1)}\cdot\dot{\boldsymbol{v}}_{\xi\xi}) \\
&\quad\;
 - \frac13 (f^{(1)})^{-1}f_\xi^{(1)}
 \boldsymbol{v}^{(1)}\times(\boldsymbol{v}^{(1)}\times\dot{\boldsymbol{v}}_{\xi\xi}).
\end{align*}
This together with the identity \( \dot{\boldsymbol{z}}_\xi
 = \boldsymbol{v}^{(1)}\times\dot{\boldsymbol{v}}_{\xi\xi} + \boldsymbol{v}_\xi^{(1)}\times\dot{\boldsymbol{v}}_\xi \) 
yield 
\begin{equation}\label{eqdz}
\dot{\boldsymbol{z}}_t
= f^{(1)}\boldsymbol{v}^{(1)}\times\bigl( \dot{\boldsymbol{z}}_{\xi\xi}
 - 2(f^{(1)})^{2/3}(\boldsymbol{v}^{(1)}\times\boldsymbol{v}_\xi^{(1)}) )
  (\boldsymbol{v}^{(1)}\cdot\dot{\boldsymbol{v}}_{\xi\xi}) \bigr)
 + \frac73 f_\xi^{(1)}\boldsymbol{v}^{(1)}\times\dot{\boldsymbol{z}}_\xi + \dot{\boldsymbol{G}}_2,
\end{equation}
where 
\begin{align*}
\dot{\boldsymbol{G}}_2
&= \boldsymbol{v}^{(1)}\times( \dot{\boldsymbol{G}}
 - f^{(1)}\boldsymbol{v}_{\xi\xi}^{(1)}\times\dot{\boldsymbol{v}}_\xi
 - 2f_\xi^{(1)}\boldsymbol{v}_\xi^{(1)}\times\dot{\boldsymbol{v}}_\xi )
 + \boldsymbol{v}_t^{(1)}\times\dot{\boldsymbol{v}}_\xi
 - \frac13 f_\xi^{(1)}\boldsymbol{v}^{(1)}\times(\boldsymbol{v}_\xi^{(1)}\times\dot{\boldsymbol{v}}_\xi),
\end{align*}
which satisfies 
\( \|\dot{\boldsymbol{G}}_2\|_1 
 \leq C( \|\dot{\boldsymbol{v}}\|_2 + \|\dot{\boldsymbol{x}}\| )
 \leq C\|\dot{\boldsymbol{x}}\|_3 \). 
In view of \eqref{eqdz} we introduce a new variable $\dot{\boldsymbol{u}}$ by 
\begin{equation}
\dot{\boldsymbol{u}} := \dot{\boldsymbol{z}}
 - 2(f^{(1)})^{2/3}(\boldsymbol{v}^{(1)}\times\boldsymbol{v}_\xi^{(1)})
   (\boldsymbol{v}^{(1)}\cdot\dot{\boldsymbol{v}}). 
\end{equation}
Then, we have 
\begin{equation}\label{eqdu}
\dot{\boldsymbol{u}}_t = f^{(1)}\boldsymbol{v}^{(1)}\times\dot{\boldsymbol{u}}_{\xi\xi}
 + \frac73 f_\xi^{(1)}\boldsymbol{v}^{(1)}\times\dot{\boldsymbol{u}}_\xi + \dot{\boldsymbol{G}}_3,
\end{equation}
where 
\begin{align*}
\dot{\boldsymbol{G}}_3
&= - 2(f^{(1)})^{2/3}(\boldsymbol{v}^{(1)}\times\boldsymbol{v}_\xi^{(1)})
   (\boldsymbol{v}^{(1)}\cdot\dot{\boldsymbol{v}}_t) 
 - 2( (f^{(1)})^{2/3}(\boldsymbol{v}^{(1)}\times\boldsymbol{v}_\xi^{(1)})
  \otimes\boldsymbol{v}^{(1)} )_t \dot{\boldsymbol{v}} \\
&\quad\;
 + 2f^{(1)}\boldsymbol{v}^{(1)}\times\bigl( 
  [\partial_\xi^2, (f^{(1)})^{2/3}(\boldsymbol{v}^{(1)}\times\boldsymbol{v}_\xi^{(1)})\otimes\boldsymbol{v}^{(1)}]
  \dot{\boldsymbol{v}} \bigr) \\
&\quad\;
 + \frac{14}{3}f_\xi^{(1)}\boldsymbol{v}^{(1)}\times\bigl(
  (f^{(1)})^{2/3}(\boldsymbol{v}^{(1)}\times\boldsymbol{v}_\xi^{(1)})
   (\boldsymbol{v}^{(1)}\cdot\dot{\boldsymbol{v}}) \bigr)_\xi + \dot{\boldsymbol{G}}_2.
\end{align*}
Here, we see that 
\begin{align*}
\boldsymbol{v}^{(1)}\cdot\dot{\boldsymbol{v}}_t
&= \boldsymbol{v}^{(1)}\cdot(\boldsymbol{v}_t^{(1)} - \boldsymbol{v}_t^{(2)}) \\
&= \boldsymbol{v}^{(1)}\cdot\bigl\{
 (f^{(1)}\boldsymbol{v}^{(1)}\times\boldsymbol{v}_{\xi\xi}^{(1)}
 + f_{\xi}^{(1)}\boldsymbol{v}^{(1)}\times\boldsymbol{v}_{\xi}^{(1)}
 + {\rm D}\boldsymbol{F}(\boldsymbol{x}^{(1)},t)\boldsymbol{v}^{(1)} ) \\
&\qquad
 - (f^{(2)}\boldsymbol{v}^{(2)}\times\boldsymbol{v}_{\xi\xi}^{(2)}
 + f_{\xi}^{(2)}\boldsymbol{v}^{(2)}\times\boldsymbol{v}_{\xi}^{(2)}
 + {\rm D}\boldsymbol{F}(\boldsymbol{x}^{(2)},t)\boldsymbol{v}^{(2)} ) \bigr\} \\
&= - \dot{\boldsymbol{v}}\cdot( 
 f^{(2)}\boldsymbol{v}^{(2)}\times\boldsymbol{v}_{\xi\xi}^{(2)}
 + f_{\xi}^{(2)}\boldsymbol{v}^{(2)}\times\boldsymbol{v}_{\xi}^{(2)} ) \\
&\quad\;
 + \boldsymbol{v}^{(1)}\cdot\bigl(
 {\rm D}\boldsymbol{F}(\boldsymbol{x}^{(1)},t)\boldsymbol{v}^{(1)} 
 - {\rm D}\boldsymbol{F}(\boldsymbol{x}^{(2)},t)\boldsymbol{v}^{(2)} \bigr),
\end{align*}
so that \( \dot{\boldsymbol{G}}_3 \) satisfies the estimate 
\( \|\dot{\boldsymbol{G}}_3\|_1 
 \leq C( \|\dot{\boldsymbol{v}}\|_2 + \|\dot{\boldsymbol{x}}\| )
 \leq C\|\dot{\boldsymbol{x}}\|_3 \). 
We further introduce a new variable \( \dot{\boldsymbol{w}} \) by 
\( \dot{\boldsymbol{u}} = a^{(1)}\dot{\boldsymbol{w}} \) with a scalar function $a^{(1)}$ defined by 
\( a^{(1)} = (f^{(1)})^{-7/6} \). 
Then, we have 
\begin{equation}\label{eqdw}
\dot{\boldsymbol{w}}_t = f^{(1)}\boldsymbol{v}^{(1)}\times\dot{\boldsymbol{w}}_{\xi\xi}
 + \dot{\boldsymbol{G}}_4,
\end{equation}
where 
\[
\dot{\boldsymbol{G}}_4 = (a^{(1)})^{-1}\biggl\{ 
 \biggl( f^{(1)}a_{\xi\xi}^{(1)} + \frac73 f_\xi^{(1)} a_\xi^{(1)} \biggr)
  \boldsymbol{v}^{(1)}\times\dot{\boldsymbol{w}}
 + \dot{\boldsymbol{G}}_3 - a_{t}^{(1)}\dot{\boldsymbol{w}} \biggr\},
\]
which satisfies \( \|\dot{\boldsymbol{G}}_4\|_1 \leq C\|\dot{\boldsymbol{x}}\|_3 \). 
Therefore, we obtain 
\begin{equation}\label{estdw}
\frac{\rm d}{{\rm d}t} \|\dot{\boldsymbol{w}}_\xi\|^2 \leq C\|\dot{\boldsymbol{x}}\|_3^2.
\end{equation}

On the other hand, we can directly obtain 
\begin{equation}\label{estdx}
\frac{\rm d}{{\rm d}t} \|\dot{\boldsymbol{x}}\|^2 \leq C\|\dot{\boldsymbol{x}}\|_2^2.
\end{equation}
Now, we define \( \dot{E}(t) \) by 
\[
\dot{E}(t) = \| \dot{h}_\xi+\boldsymbol{v}_\xi^{(1)}\cdot\dot{\boldsymbol{v}}_\xi \|^2
 + \|\dot{\boldsymbol{w}}_\xi\|^2 + \|\dot{\boldsymbol{x}}\|^2,
\]
which is equivalent to \( \|\dot{\boldsymbol{x}}\|_3 \). 
Adding \eqref{estdh}, \eqref{estdw}, and \eqref{estdx}, we have 
\( \frac{\rm d}{{\rm d}t} \dot{E}(t) \leq C \dot{E}(t) \), 
so that Gronwall's inequality yields \( \dot{\boldsymbol{x}} = \boldsymbol{0} \), that is, 
\( \boldsymbol{x}^{(1)} = \boldsymbol{x}^{(2)} \). 
The proof of the uniqueness of the solution is complete. 


\section{Conclusions and Discussions}
\setcounter{equation}{0}

We proved the time-local solvability of the initial value problem (\ref{ext}) in the Sobolev space 
\( H^{m}(\mathbf{T})\) for \( m\geq 5\) together with the uniqueness of the solution.

Since we didn't assume any structural conditions on 
the external flow \( \mbox{\mathversion{bold}$F$}\),
our time-local existence theorem 
has the potential to be utilized in
the mathematical analysis of the motion of vortex filaments 
in various physical situations. One such example which the
authors would like to consider is the motion of 
a pair of interacting vortex filaments. By regarding the effect of 
the induced flow of one filament on the other filament as an 
external flow, we can formulate the problem of the 
interaction of two filaments within the framework of 
the problem considered in this paper.

\section*{Acknowledgements}
The authors would like to thank the editors and 
referees for their fruitful advices to the first version 
of this paper.

\vspace*{1cm}

\noindent
Masashi Aiki\\
Department of Mathematics\\
Faculty of Science and Technology, Tokyo University of Science\\
2641 Yamazaki, Noda, Chiba 278-8510, Japan\\
E-mail: aiki\verb|_|masashi\verb|@|ma.noda.tus.ac.jp

\vspace*{5mm}
\noindent
Tatsuo Iguchi\\
Department of Mathematics\\
Faculty of Science and Technology, Keio University\\
3-14-1 Hiyoshi, Kohoku-ku, Yokohama, 223-8522, Japan\\
E-mail: iguchi\verb|@|math.keio.ac.jp

\end{document}